\pdfoutput=1
\documentclass[12pt,amsmath,amssymb,longbibliography,nofootinbib,showkeys,eprint]{revtex4-1}
\usepackage[activate={true,nocompatibility},final,tracking=true,kerning=true,spacing=true]{microtype}
\usepackage{filecontents} 
\usepackage{epigraph}
\usepackage{graphicx}
\usepackage{epstopdf}
\usepackage{braket}
\usepackage{bm}
\usepackage{numprint}
\usepackage{float}
\usepackage[T1,T2A]{fontenc}
\usepackage[utf8]{inputenc}
\usepackage{xr}
\usepackage{amsthm}
\theoremstyle{definitions}

\usepackage{calc}
\usepackage{seqsplit}
\usepackage{hyperref}

\begin{document}
\preprint{V.M.}
\title{On Numerical Estimation of
  Joint Probability Distribution from Lebesgue Integral Quadratures
}
\author{Vladislav Gennadievich \surname{Malyshkin}} 
\email{malyshki@ton.ioffe.ru}
\affiliation{Ioffe Institute, Politekhnicheskaya 26, St Petersburg, 194021, Russia}

\date{July, 19, 2018}

\begin{abstract}
\begin{verbatim}
$Id: JointDistributionLebesgueQuadratures.tex,v 1.55 2020/11/30 15:49:52 mal Exp $
\end{verbatim}
An important application of Lebesgue integral quadrature\cite{ArxivMalyshkinLebesgue} is developed.
Given two random processes, $f(x)$ and $g(x)$,
two generalized eigenvalue problems can be formulated and solved.
In addition to obtaining two Lebesgue quadratures (for $f$ and $g$) from two eigenproblems,
the projections of $f$-- and $g$-- eigenvectors on each other allow
to build a joint distribution estimator,
the most general form of which is a density--matrix correlation.
Examples of the density--matrix correlation
can be
a value--correlation
$V_{f^{[i]};g^{[j]}}$, similar to a regular correlation concept,
and a new one, a probability--correlation $P_{f^{[i]};g^{[j]}}$.
If Christoffel function average is used instead of regular average
the approach can be extended 
to an estimation of joint probability of three and more random processes.
The theory is implemented numerically;
the 
\href{http://www.ioffe.ru/LNEPS/malyshkin/code_polynomials_quadratures.zip}{software}
is available under the GPLv3 license.
\end{abstract}
\maketitle

\section{\label{intro}Introduction}
The Gaussian quadrature relatively a measure can be viewed
as an optimal discrete interpolation
of the measure; an application
of the Gaussian quadrature to a function
can be viewed as a Riemann integral sum.
In \cite{ArxivMalyshkinLebesgue} a new type of quadrature,
producing the Lebesgue integral was introduced and applied
to the problem of optimal discretization of a random process.
In this paper a new application of the Lebesgue quadrature is developed.
Assume one has two or more
random processes $f(x)$ and $g(x)$.
An application of
the Lebesgue quadrature to each of them
gives a set of (eigenvalue, eigenvector) pairs for $f$ and $g$.
Projecting $f$-- and $g$-- eigenvectors on each other
allows an estimator of joint $(f,g)$ distribution to be obtained:
$f$ and $g$ to have $n$ levels each
(equal to their Lebesgue quadrature value--nodes $f^{[i]}$ and $g^{[i]}$),
the projections of eigenvectors define joint probability
$\left(f= f^{[i]}\right) \cap \left( g= g^{[j]}\right)$,
the most general form of which is a density--matrix correlation.
The examples of the density--matrix correlation
can be, introduced in \cite{ArxivMalyshkinMuse},
the value--correlation
$V_{f^{[i]};g^{[j]}}$, similar to a regular correlation concept,
and a new one, a probability--correlation $P_{f^{[i]};g^{[j]}}$.
A problem of three--processes joint distribution
can be approached
and
an estimation of joint probability of three and more random processes
obtained
if Christoffel function average is used instead of regular average.

\section{\label{JointDistribution} Joint distribution estimation}
In the work \cite{ArxivMalyshkinLebesgue} a concept
of the Lebesgue integral quadrature was introduced
and optimal estimator of  random process $f(x)$ distribution
was obtained
by solving generalized eigenvalue  problem\footnote{
An important feature of the approach is its applicability\cite{liionizversiyaran,malyshkin2018spikes}
to the signals with spikes, fat tails, infinite $\Braket{f^2}$, $\Braket{fg}$, or $\Braket{g^2}$, etc.,
because eigenvalue problem (\ref{GEV}) is stable and well defined
for such $f(x)$ and $g(x)$.
The approach cannot be applied
 to the processes with
  infinite $\Braket{f}$ or $\Braket{g}$ as  this leads to infinite
  $\Braket{Q_j|f|Q_k}$ or $\Braket{Q_j|g|Q_k}$.
  }:
\begin{eqnarray}
&&\sum\limits_{k=0}^{n-1} \Braket{Q_j|f|Q_k} \alpha^{f;[i]}_k =
  \lambda_f^{[i]} \sum\limits_{k=0}^{n-1} \Braket{ Q_j|Q_k} \alpha^{f;[i]}_k
  \label{GEV} \\ 
&&  \psi_{f}^{[i]}(x)=\sum\limits_{k=0}^{n-1} \alpha^{f;[i]}_k Q_k(x)
  \label{psiC}
\end{eqnarray}
The value--nodes and weights of Lebesgue integral quadrature:
\begin{eqnarray}
  f^{[i]}&=&\lambda_f^{[i]} \label{fiLeb} \\
  w^{[i]}_f&=& \Braket{\psi_{f}^{[i]}}^2 \label{wiLeb}
\end{eqnarray}
is optimal $n$--point discrete distribution,
producing Lebesgue integral $df$ relatively the measure $d\mu$
(if $f(x)=x$ then $f^{[i]}$ and $w^{[i]}_f$ are the nodes and the weights
of Gaussian quadrature, optimal
\href{https://en.wikipedia.org/wiki/Riemann_sum}{Riemann integral sum}).
The weights $w^{[i]}_f$ of the
\href{https://en.wikipedia.org/wiki/Lebesgue_integration#Intuitive_interpretation}{Lebesgue integral}
quadrature
is the measure of $f= f^{[i]}$ sets.
Now assume one has two random processes $f(x)$ and $g(x)$
and want to obtain not only $f$-- and $g$-- distributions
(from the Lebesgue integral quadrature),
but also to obtain the measure of $\left(f= f^{[i]}\right) \cap \left( g= g^{[j]}\right)$ sets,
a joint distribution of $(f,g)$.

As we emphasized in \cite{ArxivMalyshkinLebesgue}
any quadrature is defined by $n$ (eigenvalue, eigenvector)
pairs of (\ref{GEV}) problem: $(\lambda_f^{[i]},\psi_{f}^{[i]})$.
One can
solve (\ref{GEV}) with $g(x)$ (instead of $f(x)$)
and obtain one more
set of $n$ pairs $(\lambda_{g}^{[i]},\psi_{g}^{[i]})$.
Because right hand side of (\ref{GEV}) is the same for $f(x)$ and $g(x)$,
the eigenvectors for $f$ and $g$
can be projected to each other as: $\Braket{\psi_{f}^{[i]}|\psi_{g}^{[j]}}=\sum_{m,m^{\prime}=0}^{n-1}\alpha^{f;[i]}_m\Braket{ Q_m|Q_{m^{\prime}}} \alpha^{g;[j]}_{m^{\prime}}$

Consider total measure $\int d\mu$,
taking into account $1=\sum_{i=0}^{n-1}\Braket{\psi^{[i]}}\psi^{[i]}(x)$ obtain:
\begin{align}
\int d\mu &=  \Braket{1}=
\sum_{i,j=0}^{n-1}\Braket{\psi_{f}^{[i]}}\Braket{\psi_{f}^{[i]}|\psi_{g}^{[j]}}\Braket{\psi_{g}^{[j]}} \label{int2} \\
V_{f^{[i]};g^{[j]}}&=\Braket{\psi_{f}^{[i]}}\Braket{\psi_{f}^{[i]}|\psi_{g}^{[j]}}\Braket{\psi_{g}^{[j]}}
\label{Vm} \\
\sum_{i,j=0}^{n-1}V_{f^{[i]};g^{[j]}}&=\Braket{1} \label{wsum1}\\
\sum_{j=0}^{n-1}V_{f^{[i]};g^{[j]}}&=w^{[i]}_f \label{wsum2}\\
\sum_{i=0}^{n-1}V_{f^{[i]};g^{[j]}}&=w^{[j]}_g \label{wsum3}
\end{align}

When $f=g$ then the matrix $V_{f^{[i]};g^{[j]}}$ is diagonal,
the diagonal elements are equal to the Lebesgue integral quadrature weights (\ref{wiLeb}).
In general case (\ref{Vm}) matrix elements can be considered as
Lebesgue measure of
$\left(f= f^{[i]}\right) \cap \left(  g= g^{[j]}\right)$ sets,
this matrix is related to \textbf{joint distribution matrix} of $(f,g)$
relatively the measure $d\mu$.
However, in contrast with the Lebesgue integral quadrature weights $w^{[i]}_f$ and $w^{[i]}_g$,
the $V_{f^{[i]};g^{[j]}}$ elements are not always positive,
but are still useful for $(f,g)$ joint probability estimation;
obtained joint distribution matrix
is correct ``on  average'':  (\ref{wsum2}) and (\ref{wsum3}).

One more matrix $P_{ij}$ can be introduced:
\begin{align}
P_{f^{[i]};g^{[j]}}&=\Braket{\psi_{f}^{[i]}|\psi_{g}^{[j]}}^2
\label{Pm} \\
\sum_{i,j=0}^{n-1}P_{f^{[i]};g^{[j]}}&=n
\label{wsumP}
\end{align}
The $P_{f^{[i]};g^{[j]}}$ matrix is always positive
but is normalized to eigenvalues number $n$, not to the total measure
$\Braket{1}$ as $V_{f^{[i]};g^{[j]}}$ (\ref{wsum1}) is.

The $V_{f^{[i]};g^{[j]}}$ and $P_{f^{[i]};g^{[j]}}$ joint distributions matrices
are a generalization of
value--correlation and probability--correlation
concepts introduced in Ref. \cite{ArxivMalyshkinMuse}
for  $n=2$ case.
For $n=2$ the $(P_{f^{[0]};g^{[0]}}+P_{f^{[1]};g^{[1]}}-P_{f^{[0]};g^{[1]}}-P_{f^{[1]};g^{[0]}})\big/2$
is exactly the $\widetilde{\rho}(f,g)$ of Ref. \cite{ArxivMalyshkinMuse} Appendix C,
$P_{f^{[i]};g^{[j]}}$ has the meaning of \textsl{probability of probability} for $f=f^{[i]}$ and $g=g^{[j]}$;
$V_{f^{[i]};g^{[j]}}$ has the meaning
of \textsl{probability} (measure) for $f=f^{[i]}$ and $g=g^{[j]}$
and normalizing (\ref{wsum1}),
it is different from $L^4covariation_{f^{[i]},g^{[j]}}$
of Ref. \cite{ArxivMalyshkinMuse} Appendix B,
only in basis choice. These two correlation concepts are special
case of the density matrix correlation concept.

\subsection{\label{jointdistribGeneral}Density Matrix Correlation}
Obtained (\ref{Vm}) and (\ref{Pm}) joint probability estimators can be generalized
as
\begin{eqnarray}
P_{f^{[i]};g^{[j]}}&=&\Braket{\psi_{f}^{[i]}|\psi_{g}^{[j]}} \Braket{\psi_{f}^{[i]}|\rho|\psi_{g}^{[j]}}
\label{PmDM} \\
\sum_{i,j=0}^{n-1}P_{f^{[i]};g^{[j]}}&=&\mathrm{Spur}\, \|\rho\| \label{wsumDM}
\end{eqnarray}
where $\|\rho\|$ is a ``density matrix'' operator,
 $f$ and $g$ averages are now
$\mathrm{Spur}\, \|\rho|f\|$ and $\mathrm{Spur}\, \|\rho|g\|$,
with normalizing  (\ref{wsumDM}).
The $\|\rho\|=\Ket{1}\Bra{1}$ gives regular average and value--correlation (\ref{Vm}),
the $\|\rho\|=\|1\|$ gives the number of eigenvalues
as average and probability--correlation
(\ref{Pm}).
For a given polynomial $P(x)$ of $2n-2$ degree,
density matrix $\|\rho\|$ operator can be
constructed according to Ref.\cite{ArxivMalyshkinLebesgue} Appendix A algorithm,
where, for a given polynomial, an operator $\|\rho\|$  has been obtained, such that:
\begin{eqnarray}
  \rho(x,y)&=&\sum_{i=0}^{n-1}\psi^{[i]}(x) \lambda^{[i]} \psi^{[i]}(y) \label{rhoP} \\
  P(x)&=& \rho(x,x) \label{rhoPrelation}
\end{eqnarray}
with $\lambda^{[i]} ; \psi^{[i]}(x)$
the eigenvalues and the eigenvectors of $\|\rho\|$.
In a general case the sign of (\ref{PmDM}) matrix elements, same as for (\ref{Vm}),
is not always positive what allows only ``on average'' interpretation.
One may consider a different type of estimation:
\begin{eqnarray}
P_{f^{[i]};g^{[j]}}&=&\Braket{\psi_{f}^{[i]}|\rho|\psi_{g}^{[j]}}^2
\label{PmDM2}
\end{eqnarray}
where the density matrix is used in both terms in (\ref{PmDM}).
For $\|\rho\|=\Ket{1}\Bra{1}$
the weight of $\left(f= f^{[i]}\right) \cap \left( g= g^{[j]}\right)$
is equal to the product of Lebesgue quadrature weights for $f$ and $g$ quadratures:
\begin{eqnarray}
  P_{f^{[i]};g^{[j]}}&=&\Braket{\psi_{f}^{[i]}}^2\Braket{\psi_{g}^{[j]}}^2 \label{fguncorr}\\
  \sum_{i,j=0}^{n-1}P_{f^{[i]};g^{[j]}}&=&\Braket{1}^2 \label{mu2w}
\end{eqnarray}
with (\ref{mu2w}) normalizing. The (\ref{fguncorr}) is ``uncorrelated'' answer,
where the probability of joint $(f,g)$
distribution  is equal to the product of individual distributions.
Such ``uncorrelated'' answers always arise for pure states, the states
with the density matrix in $\|\rho\|=\Ket{\psi}\Bra{\psi}$ form.
In general case, such as (\ref{rhoP}),
obtained joint distribution does not factorize.

\section{\label{numerical}Numerical Estimation}
Numerical estimation is not much different from 
the Lebesgue integral quadrature calculation in Ref. \cite{ArxivMalyshkinLebesgue}.
With a good choice for $Q_k(x)$ basis the problem can
be efficiently solved\cite{2015arXiv151005510G,2015arXiv151101887G} for a very large $n$.
Once the $\Braket{Q_m}$, $\Braket{fQ_m}$, and $\Braket{gQ_m}$
moments are obtained for $m=0\dots 2n-1$, the matrices
$\Braket{Q_j|Q_k}$, $\Braket{Q_j|f|Q_k}$, and $\Braket{Q_j|g|Q_k}$, $j,k=0\dots n-1$,
can be calculated using basis functions multiplication operator,
and eigenvalue problem (\ref{GEV})
can be solved using e.g. \href{http://www.netlib.org/lapack/lug/node54.html}{generalized eigenvalue}
subroutines from Lapack\cite{lapack}.
In practical application it is better to solve generalized eigenvalue  problem
for $f(x)$ first, obtain $\psi_{f}^{[i]}(x)$. Then generalized eigenvalue  problem
for $g(x)$ can be written in $\psi_{f}^{[i]}(x)$ basis
$\sum_{k=0}^{n-1}\Braket{\psi_{f}^{[j]}|g|\psi_{f}^{[k]}}\alpha^{g;[i]}_k=
\lambda_{g}^{[i]} \sum_{k=0}^{n-1}\Braket{\psi_{f}^{[j]}|\psi_{f}^{[k]}}\alpha^{g;[i]}_k$,
that
has unit right hand side matrix:
$\Braket{\psi_{f}^{[j]}|\psi_{f}^{[k]}}=\delta_{jk}$,
and can be solved using
\href{http://www.netlib.org/lapack/lug/node48.html}{regular eigenvalue} subroutines,
what both optimizes the code and increases numerical stability of calculations.
Provided java--implementation\cite{polynomialcode},
file
\href{http://www.ioffe.ru/LNEPS/malyshkin/code_polynomials_quadratures.zip}{\texttt{\seqsplit{code\_polynomials\_quadratures.zip}}},
the methods \texttt{\seqsplit{getValueCorrelation}},
\texttt{\seqsplit{getProbabilityCorrelation}},
and
\texttt{\seqsplit{getDensityMatrixCorrelation}}
of \texttt{\seqsplit{com/polytechnik/utils/LebesgueQuadratureJointDistribution.java}} class, constructed from two Lebesgue quadratures as input,
calculate correlation matrices.
Octave file \texttt{\seqsplit{LebesgueQuadratures\_call\_java\_example.m}}
is usage demonstration of this  
java implementation
application
to sampled data.
Given $l=1\dots M$ sample 
\begin{align}
  x^{(l)}&\to f^{(l)},g^{(l)} & \text{weight $\omega^{(l)}$}  \label{mlproblem2}
\end{align}
the matrices
\begin{align}
\Braket{Q_j|Q_k}&=\sum_{l=1}^{M}Q_j(x^{(l)})Q_k(x^{(l)})\omega^{(l)} \\
\Braket{Q_j|f|Q_k}&=\sum_{l=1}^{M}Q_j(x^{(l)})Q_k(x^{(l)})f^{(l)}\omega^{(l)} \\
\Braket{Q_j|g|Q_k}&=\sum_{l=1}^{M}Q_j(x^{(l)})Q_k(x^{(l)})g^{(l)}\omega^{(l)}
\end{align}
are calculated from sampled data.
For implementation demonstration
a subset of java library functionality
is also implemented in
native octave language in \texttt{\seqsplit{com/polytechnik/utils/LebesgueQuadraturesJointDistribution.m}}, see usage demonstration  in
\texttt{\seqsplit{LebesgueQuadratures.m}}.
The results of java and octave implementations are identical.

\section{\label{threeprocChristoffel}Christoffel Function Average And
 Joint Distribution of Three Random Processes.}
Obtained joint probability estimator (\ref{PmDM})
uses projections of $f$-- and $g$-- eigenvectors on each other.
This approach, besides ``negative probability'' interpretation difficulty
cannot directly estimate a joint distribution of three or more processes\footnote{
  Here we also consider $\mathbf{x}$--space to be a vector space of arbitrary origin;
the consideration of above corresponds to these components being polynomials of
some real variable: $x_k=Q_k(x)$}:
\begin{align}
  (x_0,x_1,\dots,x_k,\dots,x_{n-1})^{(l)}
  &\to f^{(l)},g^{(l)},e^{(l)},\dots & \text{weight $\omega^{(l)}$}  \label{mlproblem}
\end{align}
When eigenvalue/eigenvectors pairs $(\lambda_{s}^{[i]},\psi_{s}^{[i]})$ are obtained
for each $s=f,g,e, \dots$:
\begin{align}
  \Ket{s|\psi^{[i]}_{s}}&=\lambda_{s}^{[i]}\Ket{\psi^{[i]}_{s}}
  \label{pisifeg}
\end{align}
a three--function projection cannot be directly obtained
for say $\Ket{\psi^{[i]}_{f}}$, $\Ket{\psi^{[j]}_{g}}$, and $\Ket{\psi^{[k]}_{e}}$.
However,
a  three--function joint distribution estimator can be obtained with
a Christoffel function average;
the concept was originally introduced in \cite{malyshkin2015norm}
(as Eq. (20) therein)
and recently was applied in \cite{malyshkin2019radonnikodym}
to the Low Rank Representation (LRR) problem of a matrix.

Consider a  $\mathbf{x}=\mathbf{y}$ localized
state $\psi_{\mathbf{y}}(\mathbf{x})$:
\begin{align}
       \psi_{\mathbf{y}}(\mathbf{x})&=\frac{\sum\limits_{i=0}^{n-1}\psi^{[i]}(\mathbf{y})\psi^{[i]}(\mathbf{x})}
           {\sqrt{\sum\limits_{i=0}^{n-1}\left[\psi^{[i]}(\mathbf{y})\right]^2}}
           =
           \frac{\sum\limits_{j,k=0}^{n-1}y_jG^{-1}_{jk}x_k}
           {\sqrt{\sum\limits_{j,k=0}^{n-1}y_jG^{-1}_{jk}y_k}}             
  \label{psiYlocalized}
\end{align}
Here $G_{jk}=\Braket{x_j|x_k}$ is Gram matrix; in one--dimensional case it is
$G_{jk}=\Braket{Q_j|Q_k}$ in right hand side of (\ref{GEV}).
Averaging with $\psi^2_{\mathbf{y}}$  gives Radon--Nikodym approximation at $\mathbf{y}$.

The Christoffel function average of
a $\mathbf{x}$--dependent function $\phi(\mathbf{x})$ is defined as:
\begin{itemize}
\item Calculate the $\Braket{\psi_{\mathbf{y}}|\phi|\psi_{\mathbf{y}}}$
  according to
  the measure (\ref{mlproblem}), obtain $\mathbf{y}$--dependent function.
  This is just regular $\psi^2$ average we used here and in all our previous works.
  This average of  $\phi(\mathbf{x})$ with  $\psi^2_{\mathbf{y}}(\mathbf{x})$
  is it's Radon--Nikodym approximation at  $\mathbf{y}$.
\item Average obtained $\mathbf{y}$--dependent function 
over all $\mathbf{y}\in \mathbf{x}^{(l)}$, $l=1\dots M$
with the weights $\omega^{(l)}$:
\end{itemize}
\begin{align}
  \Braket{\phi}_{{\rm Christoffel}}&=
  \sum\limits_{l=1}^{M}\omega^{(l)} \Braket{\psi_{\mathbf{y}^{(l)}}|\phi|\psi_{\mathbf{y}^{(l)}}}=
    \sum\limits_{l=1}^{M}\omega^{(l)} \sum\limits_{l^{\prime}=1}^{M} \omega^{(l^{\prime})} \phi(\mathbf{x}^{(l^{\prime})})
    \psi^2_{\mathbf{y}^{(l)}}(\mathbf{x}^{(l^{\prime})})
    \label{Christoffelmeasure}
\end{align}
While regular averaging
\begin{align}
  \Braket{\phi}&=\sum\limits_{l=1}^{M}\omega^{(l)} \phi(\mathbf{y}^{(l)})
  \label{regularmeasure}
\end{align}
uses the value $\phi(\mathbf{y}^{(l)})$ to average, the Christoffel function
average (\ref{Christoffelmeasure}) uses it's
Radon--Nikodym approximation $\Braket{\psi_{\mathbf{y}^{(l)}}|\phi|\psi_{\mathbf{y}^{(l)}}}$ instead.
For a large enough $n$ the result is similar.
Since Radon--Nikodym approximation preserves the normalizing and sign
we also have: $\Braket{1}_{{\rm Christoffel}}=\Braket{1}$.

The key Christoffel average feature is that we can represent it in any basis: because the $\psi_{\mathbf{y}}(\mathbf{x})$ is a regular
wavefunction, it can be expanded in any full basis, for example (\ref{pisifeg}).
For  $\psi^{[i]}_{s}(\mathbf{x})$,  $s=f,g,e, \dots$, taking into account basis functions
orthogonality, obtain:
\begin{align}
  \psi_{\mathbf{y}}(\mathbf{x})&=  \sum\limits_{i=0}^{n-1}\Braket{\psi^{[i]}_{s}|\psi_{\mathbf{y}}}\psi^{[i]}_{s}(\mathbf{x})
  \label{psilocalizedexpansion} \\  
  \Braket{s}_{{\rm Christoffel}}&=
  \sum\limits_{l=1}^{M}\omega^{(l)} \Braket{\psi_{\mathbf{y}^{(l)}}|s|\psi_{\mathbf{y}^{(l)}}}=
  \sum\limits_{l=1}^{M}\omega^{(l)} \sum\limits_{i=0}^{n-1}\lambda_{s}^{[i]}\Braket{\psi_{\mathbf{y}^{(l)}}|\psi^{[i]}_{s}}^2
  \label{pissexmapnsion} \\
  1&= \sum\limits_{i=0}^{n-1} \Braket{\psi_{\mathbf{y}^{(l)}}|\psi^{[i]}_{s}}^2
  \label{normalize}
\end{align}
The Christoffel function average $\Braket{s}_{{\rm Christoffel}}$
is an average of $l=1\dots M$ observations,
where each one is a probability distribution of $n$ outcomes $\lambda_{s}^{[i]}$ (the (\ref{pisifeg}) eigenvalues) with the weights $\omega^{(l)}\Braket{\psi_{\mathbf{y}^{(l)}}|\psi^{[i]}_{s}}^2$; the regular average $\Braket{s}$ is a superposition of $l=1\dots M$ observations $s^{(l)}$ with the weights $\omega^{(l)}$.

A remarkable feature of (\ref{pissexmapnsion}) is that it has
the same representation for any $s=f,g,e,\dots$.
The most straightforward approach to obtain joint distribution 
is, for a given $l$, to consider the components to be independent
(this is reasonable because for a large $n$
only a single coefficient
$\Braket{\psi_{\mathbf{y}^{(l)}}|\psi^{[k]}_{s}}^2$ is typically large)
hence
one can consider the joint probability 
at $\mathbf{y}^{(l)}$ to be a product\footnote{
  One can also try to consider the ``two--particle'' probability distribution
  like (\ref{Pm}) with
  \href{https://en.wikipedia.org/wiki/Chain_rule_(probability)}{chain rule}
  and decoupling.}
$\Braket{\psi_{\mathbf{y}^{(l)}}|\psi^{[i]}_{f}}^2\Braket{\psi_{\mathbf{y}^{(l)}}|\psi^{[j]}_{g}}^2
\Braket{\psi_{\mathbf{y}^{(l)}}|\psi^{[k]}_{e}}^2$. The joint
probability is then:
\begin{align}
  P_{f^{[i]};g^{[j]};e^{[k]}} &=
  P \left(f= f^{[i]} \cap  g= g^{[j]} \cap  e= e^{[k]} \right)
  =\sum\limits_{l=1}^{M}\omega^{(l)}
  \Braket{\psi_{\mathbf{y}^{(l)}}|\psi^{[i]}_{f}}^2
  \Braket{\psi_{\mathbf{y}^{(l)}}|\psi^{[j]}_{g}}^2
  \Braket{\psi_{\mathbf{y}^{(l)}}|\psi^{[k]}_{e}}^2
  \label{Prob3V}
\end{align}
This joint probability estimator is always positive,
has proper normalizing,
and can be used to estimate joint distribution of three or more random processes.
In one--dimensional case it gives outcome weights
equal to Christoffel function diagonal elements
in the basis $\Ket{\psi^{[i]}_{f}}$.
For $f=g=e$ it produces properly diagonal
joint distribution
matrix
 $P \left(f= f^{[i]} \cap  g= g^{[j]} \cap  e= e^{[k]} \right)\approx 0$ for $i\ne j \ne k$
only for a large enough $n$, this is a limitation of this estimator.

We see (\ref{Prob3V}) as an extremely promising
path to vector--valued class label machine learning.
The Chrisoffel function joint probability estimator
combines the spectral approach (\ref{pisifeg})
success for scalar class label of Ref. \cite{malyshkin2019radonnikodym}
with a vector class label $\mathbf{f}=(f,g,e,\dots)$.

\section{\label{conclusion}Conclusion}
Obtained in Ref.\cite{ArxivMalyshkinLebesgue} a new class of quadratures,
the Lebesgue quadrature,
can be applied not only to optimal discretization
of a random process by a $n$--point discrete Lebesgue measure,
but also to a numerical estimation of joint distribution of $(f(x),g(x))$.
The most general form is
density matrix correlation.
Introduced in 
Ref. \cite{ArxivMalyshkinMuse} Appendix B,
value--correlation (\ref{Vm})
and probability--correlation (\ref{Pm}) are
special cases of density matrix correlation.
If Christoffel function average is used instead of regular average
the approach can be further extended (\ref{Prob3V})
to three and more random processes.
\href{http://www.ioffe.ru/LNEPS/malyshkin/code_polynomials_quadratures.zip}{The software}
 is available under the GPLv3 license.

\bibliography{LD}

\end{document}